\documentclass[leqno,12pt]{amsart}

\usepackage{amsmath,amstext,amssymb,amsopn,amsthm,mathrsfs}
\usepackage{bbm}
\usepackage{verbatim}
\usepackage{enumitem}
\usepackage{color}

\newcommand*\RR{\mathbb{R}}
\newcommand*\NN{\mathbb{N}}
\newcommand*\ind{\mathbbm{1}}
\newcommand*\ven{\vert n\vert}
\newcommand*\al{\alpha}
\newcommand*\fun{\varphi_n}

\newcommand*\hfun{\varphi_n^\al}

\newcommand*\hfunk{\varphi^{\alpha}_{k}}
\newcommand*\Rop{R_r^{\alpha}}

\newcommand*\ve{\varepsilon}
\newcommand*\mf{\mathfrak{m}}

\title[On Hardy's inequality for Hermite expansions]
{On Hardy's inequality for Hermite expansions}

\author[P{.} Plewa]{Pawe\l{} Plewa}
\address{Pawe\l{} Plewa \newline
			Faculty of Pure and Applied Mathematics, 
      Wroc\l{}aw University of Science and Technology       \newline
      Wyb{.} Wyspia\'nskiego 27,
      50--370 Wroc\l{}aw, Poland      
      }
\email{pawel.plewa@pwr.edu.pl}

\allowdisplaybreaks

\theoremstyle{plain}
\newtheorem{thm}{Theorem}[section]

\newtheorem{lm}[thm]{Lemma}
\newtheorem{prop}[thm]{Proposition}

\theoremstyle{definition}

\theoremstyle{remark}
\newtheorem*{rem*}{Remark}

\theoremstyle{plain}

\usepackage[dvipsnames]{xcolor}

\newcounter{comcount}

\begin{document}
\begin{abstract}
Sharp multi-dimensional Hardy's inequality for the Laguerre functions of Hermite type is proved for the type parameter $\al\in[-1/2,\infty)^d$. As a consequence we obtain the corresponding result for the generalized Hermite expansions. In particular, it validate that the known version of Hardy's inequality for the Hermite functions is sharp.
\end{abstract}

\maketitle
\footnotetext{
The paper is a part of author's doctoral thesis written under the supervision of Professor Krzysztof Stempak.
}

\section{Introduction}
Hardy and Littlewood \cite{HardyLittlewood} proved the following inequality for Fourier coefficients
\begin{equation}\label{old_Hardy_inequality}
\sum_{k\in\mathbb{Z}}\frac{\vert \hat{f}(k)\vert}{\vert k\vert+1}\lesssim \Vert f\Vert_{\mathrm{ Re} H^1}, 
\end{equation}
where $\mathrm{ Re} H^1$ denotes the real Hardy space constituted by the boundary values of the real parts of functions in the Hardy space $H^1(\mathbb{D})$, where $\mathbb{D}$ is the unit disk on the plane.

Kanjin \cite{Kanjin1} initiated investigation of analogues of \eqref{old_Hardy_inequality} for orthogonal expansions. He proved the one-dimensional version of the following inequality
\begin{equation}\label{Hardy_Hermite_introd}
\sum_{n\in\mathbb{N}^d}\frac{\vert \langle f,h_n\rangle\vert}{ (n_1+\ldots+n_d+1)^{E}}\lesssim \Vert f\Vert_{H^1(\RR^d)},\qquad f\in H^1(\RR^d),
\end{equation}
where $n=(n_1,\ldots,n_d)$, $\langle \cdot,\cdot\rangle$ stands for the inner product in $L^2(\RR^d)$, $\{h_n\}_{n\in\NN^d}$ are the Hermite functions, and $H^1(\RR^d)$ denotes the Hardy space. We will refer to the constant $E$ as the admissible exponent. 

Recently many authors studied Hardy's inequality for Hermite expansions. In the mentioned article Kanjin examined only the case $d=1$ and proved a version of \eqref{Hardy_Hermite_introd} with $E=29/36$. Later Radha \cite{Radha} investigated the multi-dimensional setting $d\geq 1$. For an arbitrary $\ve>0$, the admissible exponent $E=(17d+12+\ve)/(12d+24)$ was obtained. Then Radha and Thangavelu \cite{RadhaThangavelu} received $E=3d/4$ for $d\geq 2$. Unfortunately, the applied method did not work in the one-dimensional case. Kanjin \cite{Kanjin2} basing on a paper of Balasubramanian and Radha \cite{BalasRadha} justified that for $d=1$ the admissible exponent is $E=3/4+\ve$, for an arbitrary $\ve>0$. He also conjectured that it can be lowered to $3/4$. It was indeed proved by  Z. Li, Y. Yu and Y. Shi \cite{LiYuShi}.

Hardy's inequality was also investigated in the context of different orthonormal expansions as well. Kanjin and Sato \cite{KanjinSato} studied the case of the Jacobi expansions. Moreover, the author considered various Laguerre expansions in \cite{Plewa, Plewa2}. Furthermore, an analogue of Hardy's inequality was scrutinized, namely the Hardy space $H^1$ was replaced by $H^p$ for $p\in(0,1)$ (see \cite{BalasRadha,RadhaThangavelu,Satake}).

The primary goal of this article is to prove that the admissible exponent in \eqref{Hardy_Hermite_introd} cannot be lowered. For this purpose we extend the result from \cite{Plewa} for Laguerre expansions of Hermite type, to a wider range of the type parameter, namely $\al\in[-1/2,\infty)^d$. We also construct an explicit counterexample to show that the associated admissible exponent $E=3d/4$ is sharp. Moreover, we are able to deduce the corresponding result for the generalized Hermite expansions along with its sharpness. Consequently, we get sharpness of \eqref{Hardy_Hermite_introd} with $E=3d/4$.

Our main tool in establishing Hardy's inequality is \cite[Theorem~2.2]{Plewa2}. The verification of the required conditions for the type parameter $\al\in(-1/2,1/2)$ is more complicated than for $\al\in\{-1/2\}\cup[1/2,\infty)$ (as it was implicitly done in \cite{Plewa}). In order to deduce Hardy's inequality for the generalized Hermite setting from the result for the Laguerre setting of Hermite type, we apply a decomposition of functions on $\RR^d$ with respect to its parity. Using the same method one can prove an $L^1$-analogue of Hardy's inequality (compare \cite{Kanjin2,Plewa,Plewa2}).

The organization of this paper is as follows. In Section 2 we state preliminaries, mainly some facts about the Hardy spaces, and recall \cite[Theorem~2.2]{Plewa2}. Section 3 is devoted to the Laguerre expansions of Hermite type. We present some auxiliary results leading to the verification of the assumptions of Theorem \ref{general_thm}. Furthermore, we construct the mentioned counterexample. In Section 4 we justify that Hardy's inequality for the generalized Hermite expansions follows from the corresponding result for the Laguerre functions of Hermite type.

%
\subsection*{Notation}
Throughout this paper we shall denote $\RR^d_+=(0,\infty)^d$ and $\NN_+=\NN\setminus\{0\}=\{1,2,\ldots\}$, where $d\geq 1$ is the dimension. We shall distinguish the one-dimensional variables from the multi-dimensional ones. Therefore, in the case $d=1$ we write $u,v$ for real variables and $k$ or $j$ for non-negative integers. On the other hand, in the case $d\geq 1$ we use $x=(x_1,\ldots,x_d),\ y=(y_1,\ldots,y_d)$ for real vectors, and $n=(n_1,\ldots,n_d)\in\NN^d$ for multi-indices. The Euclidean norm is denoted by $|x|$ and $|y|$, whereas $\ven=n_1+\ldots+n_d$ stand for the length of $n$. If a multi-index is constant, then we will use the bold font, e.g. $\mathbf{0}=(0,\ldots,0)$. The Laguerre type multi-index $\al=(\al_1,\ldots,\al_d)\in(-1,\infty)^d$ will be denoted by the same symbol in both cases $d=1$ and $d\geq 1$. It should be always clear from the context whether $\al$ refers to $d=1$ or $d\geq 1$. Similarly as before, $|\al|=\al_1+\ldots+\al_d$, stands for the length of the multi-index $\al$. Note that $|\al|$ may be negative. We will use the usual convention writing $x^\al=\prod_{i=1}^d x_i^{\al_i}$, $x\in\RR^d_+$. If a function $f$ is defined on $\RR^d$, then its restriction to $\RR^d_+$ is denoted by $f^+$.

The symbol $\lesssim$ stands for inequalities that hold with a multiplicative constant that may vary from line to line. Such constant may depend on parameters quantified beforehand, but not on the ones quantified afterwards. If $\lesssim$ and $\gtrsim$ hold simultaneously, then we will write $\simeq$.

\subsection*{Acknowledgement}
The author would like to express his gratitude to Professor Krzysztof Stempak for his remarks.

\section{Preliminaries}
A measurable function $f$ defined on $\RR^d$ is called $\eta$-symmetric for some $\eta\in\{0,1\}^d$, if $f$ is even with respect to every $i$-th coordinate such that $\eta_i=0$ and odd with respect to the remaining coordinates. We shall make use of the decomposition
\begin{equation*}
f=\sum_{\eta\in\{0,1\}^d} f_\eta,
\end{equation*}
where
\begin{equation*}
f_\eta(x)=2^{-d} \sum_{\epsilon\in\{-1,1\}^d} \epsilon^\eta f(\epsilon x).
\end{equation*}

The classical Hardy space $H^1(\RR^d)$ can be defined in many ways (see \cite{Stein}), e.g. given a Schwartz function $\psi$ such that $\int\psi\neq 0$, we say that a function $f\in L^1(\RR^d)$ belongs to $H^1(\RR^d)$ if and only if
\begin{equation}\label{H1_norm_def}
\Vert f\Vert_{H_m^1(\RR^d)}:=\big\Vert \sup_{t>0}|f\ast \psi_t|\big\Vert_{L^1(\RR^d)}<\infty,
\end{equation}
where $\psi_t(x)=t^{-d}\psi(x/t)$. The definition of $H^1(\RR^d)$ is independent of the chosen function $\psi$. The definition \eqref{H1_norm_def} is referred to as the maximal characterization of $H^1(\RR^d)$. We emphasize that 
\begin{equation*}
\Vert f\Vert_{L^1(\RR^d)}\lesssim \Vert f\Vert_{H_m^1(\RR^d)}.
\end{equation*}

A measurable function $a$ is called an $H^1(\RR^d)$-atom if it is supported in a Euclidean ball $B$ and satisfies the cancellation condition and the size condition, namely $\int a(x)\,dx=0$ and $\Vert a\Vert_{L^\infty(\RR^d)}\leq |B|^{-1}$, respectively, where $|B|$ denotes the Lebesgue measure of $B$. 

A function $f$ is in $ H^1(\RR^d)$ if and only if it admits an atomic decomposition, i.e. there exist a sequence of coefficients $\{\lambda_j\}_{j\in\NN}$ and a sequence of $H^1(\RR^d)$-atoms $\{a_j\}_{j\in\NN}$ such that
\begin{equation}\label{atomic_decomposition}
f(x)=\sum_{j=0}^\infty \lambda_j a_j(x),
\end{equation}
where the series is convergent in $H^1(\RR^d)$. Moreover,
$$\sum_{j=0}^\infty |\lambda_j|\lesssim \Vert f\Vert_{H_m^1(\RR^d)}. $$

We define
$$\Vert f\Vert_{H^1_{at}(\RR^d)}=\inf\sum_{j=0}^\infty |\lambda_j|, $$
where the the infimum is taken over all atomic decompositions of $f$. The norms $\Vert\cdot\Vert_{H^1_m(\RR^d)}$ and $\Vert\cdot\Vert_{H^1_{at}(\RR^d)}$ are equivalent. From now on, we shall use the latter and write simply $\Vert\cdot\Vert_{H^1(\RR^d)}$


We emphasise that for $f\in H^1(\RR^d)$ and every $\epsilon\in\{-1,1\}^d$ we have $\Vert f\Vert_{H^1(\RR^d)}=\Vert f(\epsilon\cdot)\Vert_{H^1(\RR^d)}$. Hence, for any $\eta\in\{0,1\}^d$ there is $f_\eta\in H^1(\RR^d)$ and $\Vert f_\eta\Vert_{H^1(\RR^d)}\leq\Vert f\Vert_{H^1(\RR^d)}$.

The following lemma holds.

{\lm\label{H1_lemma} If $\eta\in\{0,1\}^d$ and $f\in H^1(\RR^d)$ is $\eta$-symmetric, then $f\ind_{\RR^d_+}\in H^1(\RR^d)$. Moreover,
$$\Vert f\ind_{\RR^d_+}\Vert_{H^1(\RR^d)} \simeq\Vert f\Vert_{H^1(\RR^d)}. $$
}
\begin{proof}
Fix $\eta\in\{0,1\}^d$ and $\eta$-symmetric function $f\in H^1(\RR^d)$. We choose an atomic decomposition of $f$. Let
$$f(x)=\sum_{j=0}^\infty \lambda_j a_j(x), $$
where $a_j$'s are $H^1(\RR^d)$-atoms. Hence,
$$f(x)\ind_{\RR^d_+}(x)=f_\eta(x)\ind_{\RR^d_+}(x)= \sum_{j=0}^\infty \lambda_j2^{-d}\sum_{\epsilon\in\{-1,1\}^d} \epsilon^\eta a_j(\epsilon x)\ind_{\RR^d_+}(x). $$

In order to prove that $\Vert f\ind_{\RR^d_+}\Vert_{H^1(\RR^d)} \leq\Vert f\Vert_{H^1(\RR^d)}$ it suffices to justify that for any $H^1(\RR^d)$-atom $a$, the function
$$a_\eta(x)\ind_{\RR^d_+}(x)=2^{-d}\sum_{\epsilon\in\{-1,1\}^d} \epsilon^\eta a(\epsilon x)\ind_{\RR^d_+}(x) $$
is an $H^1(\RR^d)$-atom as well. Indeed, if the inferior of the support of $a$ does not intersect any of the hyperplanes $\langle e_i\rangle^{\perp}$, $i=1,\ldots,d$, then for all but one $\epsilon\in\{-1,1\}^d$, there is $a(\epsilon\cdot)\equiv 0$. For the remaining $\epsilon$ there holds $a(\epsilon\cdot)\ind_{\RR^d_+}=\epsilon^\eta a$, so $a_\eta\ind_{\RR^d_+}$ is an $H^1(\RR^d)$-atom. 

Let us now define
$$I=\big\{i\in\{1,\ldots,d\}\colon \mathrm{int}\,\mathrm{supp}\,a\cap \langle e_i\rangle^{\perp}\neq\emptyset\big\}, $$
where $\mathrm{int}$ denotes the interior of a set. Without any loss of generality we may assume that we have $I=\{1,\ldots,k\}$ for some $k\in\{1,\ldots,d\}$ . Then, for any $\epsilon_2\in\{-1,1\}^{d-k}$ the function
$$2^{-k}\Big( \sum_{\epsilon_1\in\{-1,1\}^k} (\epsilon_1,\epsilon_2)^\eta a((\epsilon_1,\epsilon_2)x)\Big)\ind_{\RR^d_+}(x) $$
is an $H^1(\RR^d)$-atom. Moreover, for all but one $\epsilon_2\in\{-1,1\}^{d-k}$ the function above vanishes identically. Therefore $a_\eta(x)$ is an $H^1(\RR^d)$-atom. Hence, $f\ind_{\RR^d_+}\in H^1(\RR^d)$ and
$$\Vert f\ind_{\RR^d_+}\Vert_{H^1(\RR^d)} \leq \Vert f\Vert_{H^1(\RR^d)}.$$

In order to justify the opposite estimate we notice that
$$\Vert f\Vert_{H^1(\RR^d)}=\Big\Vert \sum_{\epsilon\in\{-1,1\}^d}f(\epsilon\cdot) \ind_{\RR^d_+}\Big\Vert_{H^1(\RR^d)}=\Big\Vert \sum_{\epsilon\in\{-1,1\}^d}\epsilon^\eta f \ind_{\RR^d_+}\Big\Vert_{H^1(\RR^d)}\leq 2^d \Vert f\ind_{\RR^d_+}\Vert_{H^1(\RR^d)}. $$
This finishes the proof of the lemma.
\end{proof}

We define the Hardy space $H^1(\RR_+^d)$ as follows. A function $f\in L^1(\RR^d_+)$ belongs to $H^1(\RR_+^d)$ if there exists $g\in H^1(\RR^d)$ such that $\mathrm{supp}\ g\subset [0,\infty)^d$ and $g^+=f$. Moreover, we set $\Vert f\Vert_{H^1(\RR_+^d)}=\Vert g\Vert_{H^1(\RR^d)}$.

The proof of Lemma \ref{H1_lemma} yields that $f\in H^1(\RR_+^d)$ if and only if it admits an atomic decomposition as in \eqref{atomic_decomposition}, where $a_j$ are $H^1(\RR_+^d)$-atoms, e.g. $a_j$ are usual atoms and their supports are Euclidean balls intersected with $[0,\infty)^d$. Furthermore, for $\eta\in\{0,1\}^d$ and $f\in H^1(\RR^d)$ there is
\begin{equation}\label{H1_norm_estimate}
\Vert f_\eta^+\Vert_{H^1(\RR^d_+)}\simeq \Vert f_\eta\Vert_{H^1(\RR^d)}\leq \Vert f\Vert_{H^1(\RR^d)}.
\end{equation}

We shall make use of \cite[Theorem~2.2]{Plewa2}. For the reader's convenience we state it below (only for Lebesgue measure). 
{\thm\label{general_thm} Let $X$ be an open convex subset of $\RR^d$. For a given orthonormal basis $\{\fun\}_{n\in\NN^d}$ in $L^2(X)$, such that $\fun\in L^\infty(X)$, $n\in\NN^d$, we define a family of operators $\{R_r\}_{r\in(0,1)}$ via
\begin{equation}\label{R_def}
R_r f=\sum_{n\in\NN^d} r^{\ven} \langle f,\fun\rangle \fun,\qquad r\in(0,1),\qquad f\in L^2(X).
\end{equation}
We assume that the operators $R_r$ are integral operators and the associated kernels satisfy for some $\gamma>0$ and a finite set $\Delta$ composed of positive numbers the condition
\begin{equation}\label{main_condition}
\Vert R_r(x,\cdot)-R_r(x',\cdot) \Vert_{L^2(X)}\lesssim \sum_{\delta\in\Delta}|x-x'|^\delta(1-r)^{-\frac{\gamma(d+2\delta)}{d+2}},
\end{equation}
uniformly in $r\in(0,1)$, $x'\in X$, and almost every $x$ such that $|x'-x|\leq 1/3$.
Then the inequality 
$$\sum\limits_{n\in\NN^d}\frac{|\langle f,\fun\rangle|}{(\ven+1)^E}\lesssim \Vert f\Vert_{H^1(X)}, $$
holds uniformly in $f\in H^1(X)$, where
\begin{equation*}
E=\frac{\gamma d}{(d+2)}+\frac{d}{2}.
\end{equation*}
}

In the theorem above the space $H^1(X)$ is a Hardy space is the sense of Coifman-Weiss (see \cite[pp.~591-592]{CoifmanWeiss}). If $X=\RR^d$ or $X=\RR^d_+$, then it coincides with the definitions presented before.

\section{Laguerre functions of Hermite type}\label{S2}
The Laguerre functions of Hermite type are defined by the formula
$$\hfunk(u)=\Big(\frac{2\Gamma(k+1)}{\Gamma(k+\al+1)}  \Big)^{1/2}L_{k}^{\al}(u^2)u^{\al+1/2}e^{-u^2/2},\qquad u>0,$$
in the one-dimensional case, and as the tensor product in higher dimensions. The system of functions $\{\hfun\}_{n\in\NN^d}$ is an orthonormal basis in $L^2(\RR^d_+)$.

We will make use of the known estimates (see \cite[p.~435]{Muckenhoupt} and \cite[p.~699]{AskeyWainger})
\begin{equation}\label{pointwise_estimates_hfun}
\vert\hfunk(u)\vert\lesssim\left\{ \begin{array}{ll}
u^{\al+1/2}\nu^{\al/2},  & 0<u\leq 1/\sqrt{\nu},\\
\nu^{-1/4}, & 1/\sqrt{\nu}<u\leq\sqrt{\nu/2},\\
u^{1/2}(\nu(\nu^{1/3}+\vert u^2- \nu\vert))^{-1/4}, & \sqrt{\nu/2}<u\leq \sqrt{3\nu/2},\\
u^{1/2}\exp(-\gamma u^2), & \sqrt{3\nu/2}<u<\infty,
\end{array}\right.
\end{equation}
where $\nu=\nu(\al,k)=\max(4k+2\al+2,2)$ and with $\gamma>0$ depending only on $\alpha$.

Using \eqref{pointwise_estimates_hfun} for $\al\geq -1/2$ one gets
\begin{equation}\label{hfun_estim}
\Vert \hfunk\Vert_{L^\infty(\RR)}\lesssim (k+1)^{-1/12},\qquad \Vert \hfunk\Vert_{L^\infty((0,1))}\lesssim (k+1)^{-1/4},\qquad k\in\NN,
\end{equation}
compare \cite[p.~99]{StempakTohoku}. Moreover, using \eqref{pointwise_estimates_hfun} and the recurrence formula
\begin{equation}\label{derivative_formula}
\frac{d}{d u}\hfunk(u)=-2\sqrt{k}\varphi_{k-1}^{\al+1}(u)+\Big(\frac{2\al+1}{2u}-u\Big)\hfunk(u),
\end{equation}
where $\varphi^\al_{-1}\equiv 0$, we obtain for $\al\in\{-1/2\}\cup[1/2,\infty)$, 
\begin{equation}\label{hfun_der_estim}
\Big\Vert\frac{d}{d \cdot}\hfunk(\cdot)\Big\Vert_{L^\infty(\RR)}\lesssim (k+1)^{5/12},\qquad k\in\NN.
\end{equation}  
The estimate fails to hold for $\al\in(-1/2,1/2)$. However, it is easy to prove that for $\al\in[-1/2,\infty)$ we have 
\begin{equation}\label{hfun_der_estim_x>1/3}
\Big\Vert\frac{d}{d \cdot}\hfunk(\cdot)\Big\Vert_{L^\infty([1/2,\infty))}\lesssim (k+1)^{5/12},\qquad k\in\NN.
\end{equation}

In order to prove Hardy's inequality associated with the Laguerre functions of Hermite type we shall use Theorem \ref{general_thm}. The kernels associated with the family of integral operators $\{\Rop\}_{r\in(0,1)}$ for Laguerre functions of Hermite type, defined as in \eqref{R_def}, are of the form 
\begin{equation*}
\Rop(x,y)=\sum_{n\in\NN^d} r^{\vert n\vert} \hfun(x)\hfun(y),\qquad x,y\in\RR^d_+,
\end{equation*}
and, for $d=1$, can be explicitly expressed by (compare \cite[p.~102]{Szego})
\begin{equation}\label{def_R_ker_explic}
\Rop(u,v)=\frac{2 (uv)^{1/2}}{(1-r)r^{\al/2}}\exp\left(-\frac{1}{2}\frac{1+r}{1-r}(u^2+v^2)\right)I_{\al }\left(\frac{2r^{1/2}}{1-r}uv\right),
\end{equation}
where $I_{\alpha}$ denotes the modified Bessel function of the first kind, and as the tensor product in higher dimensions.

We remark that in the light of \cite[Lemma~3.1]{Plewa} in order to verify the multi-dimensional assumption \eqref{main_condition} (with $\gamma=-(d+2)/4$ and $\Delta=\{1,\al_1+1/2,\ldots,\al_d+1/2\}$) for the Laguerre functions of Hermite type with $\al\in[-1/2,\infty)^d$, it suffices to prove the following one-dimensional result.

{\prop\label{hfun_prop} If $\al\in [-1/2,\infty)$, then
$$\Big\Vert\Rop(u,\boldsymbol{\cdot})-\Rop(u',\boldsymbol{\cdot})\Big\Vert_{L^2(\RR_+)}\lesssim \frac{|u-u'|}{(1-r)^{3/4}}+\frac{|u-u'|^{\al+1/2}}{(1-r)^{(\al+1)/2}}, $$
uniformly in $r\in(0,1)$ and $u,u'>0$ such that $|u-u'|\leq 1/2$.}

Before the proof of the proposition we present two auxiliary lemmas.

{\lm\label{hfun_substr_lem} If $\al\in(-1/2,1/2)$, then
\begin{equation*}
|\hfunk(u)-\hfunk(v)|\lesssim |u-v|(k+1)^{-1/4}+|u-v|^{\al+1/2} (k+1)^{\al/2},
\end{equation*}
uniformly in $u,v\in(0,1)$ and $k\in\NN$.}
\begin{proof}
Without any loss of generality we assume that $0<u\leq v<1$. Fix $\al\in(-1/2,1/2)$ and $u,v\in(0,1)$. Note that \eqref{pointwise_estimates_hfun} yields
\begin{equation*}
\frac{|\hfunk(s)|}{s}\lesssim (k+1)^{1/4}+s^{\al-1/2}(k+1)^{\al/2},\qquad s\in(0,1),\ k\in\NN.
\end{equation*}
Hence, applying \eqref{derivative_formula}, \eqref{hfun_estim}, and  using the fact that the function $s\to s^{\al+1/2}$ is $(\al+1/2)$-H\"older continuous on $(0,1)$, we get
\begin{align*}
|\hfunk(u)-\hfunk(v)|&=\Big| \int_{u}^{v} \bigg(-2\sqrt{k}\varphi_{k-1}^{\al+1}(s)+\Big(\frac{2\al+1}{2s}-s\Big)\hfunk(s)\bigg)\,ds\Big|\\
&\lesssim |u-v|(k+1)^{1/4}+(k+1)^{\al/2}\Big|\int_{u}^{v} s^{\al-1/2}\,ds\Big|\\
&\lesssim |u-v|(k+1)^{1/4}+(k+1)^{\al/2}|u-v|^{\al+1/2},
\end{align*}
uniformly in $u,v\in(0,1)$ and $k\in\NN$. This finishes the proof.
\end{proof}

{\lm\label{hfun_lastcomponent_lemma} For $\al\in(-1/2,1/2)$ the estimate
$$\Big\Vert u^{-1}\Rop(u,\cdot)\Big\Vert_{L^2(\RR_+)}\lesssim (1-r)^{-3/4}+u^{\al-1/2}(1-r)^{-(\al+1)/2}, $$
holds uniformly in $r\in(1/2,1)$ and $u>0$.}

\begin{proof}
Fix $\al\in(-1/2,1/2)$. Using \eqref{def_R_ker_explic} and the estimates (see \cite[p.~136]{Lebedev})
$$I_\nu(s)\lesssim s^{\nu},\qquad s\in(0,1), $$
$$I_\nu(s)\lesssim s^{-1/2}e^{s},\qquad s\in(1,\infty), $$
we obtain the pointwise bound (compare \cite[(8)]{Plewa})
\begin{equation*}
\Rop(u,v)\lesssim\left\{\begin{array}{ll}
(1-r)^{-\al-1}(uv)^{\al+1/2}\exp\left(-\frac{1}{2}\frac{1+r}{1-r}(u^2+v^2) \right), & v\leq \frac{1-r}{2\sqrt{r}u},\\
(1-r)^{-1/2}\exp\left(-\frac{1}{2}\frac{1+r}{1-r}(v-u)^2 \right), & v\geq \frac{1-r}{2\sqrt{r}u}.
\end{array}\right.
\end{equation*}

Now we shall prove the claim. The following estimates are uniform in $r\in(1/2,1)$ and in the indicated ranges of $u$. Firstly, note that for $u>0$
\begin{align*}
\int_0^{\frac{1-r}{2\sqrt{r}u}}u^{-2}\Rop(u,v)^2\,dv&\lesssim (1-r)^{-2(\al+1)}\int_0^{\frac{1-r}{2\sqrt{r}u}}u^{2\al-1}v^{2\al+1}\exp\Big(-\frac{1+r}{1-r}v^2\Big)\,dv\\
&\lesssim (1-r)^{-(\al+1)}u^{2\al-1} \int_0^\infty v^{2\al+1}e^{-v^2}\,dv\\
&\lesssim (1-r)^{-(\al+1)}u^{2\al-1}.
\end{align*}
Secondly, for $u\leq (1-r)/(4\sqrt{r}u)$, we have
\begin{align*}
\int_{\frac{1-r}{2\sqrt{r}u}}^{\infty}  u^{-2}\Rop(u,v)^2\,dv&\lesssim (1-r)^{-3}\int_{\frac{1-r}{2\sqrt{r}u}}^{\infty} v^2\exp\Big(-\frac{1+r}{1-r}(v-u)^2\Big) \,dv\\
&\lesssim (1-r)^{-3} \int_{\frac{1-r}{2\sqrt{r}u}-u}^{\infty} (v+u)^2\exp\Big(-\frac{1+r}{1-r}v^2\Big) \,dv\\
&\lesssim (1-r)^{-3} \int_{0}^{\infty} v^2\exp\Big(-\frac{1+r}{1-r}v^2\Big) \,dv\\
&\lesssim (1-r)^{-3/2},
\end{align*}
and for $u\geq (1-r)/(4\sqrt{r}u)$ we obtain
\begin{align*}
\int_{\frac{1-r}{2\sqrt{r}u}}^{\infty}  u^{-2}\Rop(u,v)^2\,dv&\lesssim (1-r)^{-1}\int_{\frac{1-r}{2\sqrt{r}u}}^{\infty} (1-r)^{-1}\exp\Big(-\frac{1+r}{1-r}(v-u)^2\Big) \,dv\\
&\lesssim (1-r)^{-2} \int_{-\infty}^{\infty} \exp\Big(-\frac{1+r}{1-r}v^2\Big) \,dv\\
&\lesssim (1-r)^{-3/2}.
\end{align*}
This finishes the proof of the lemma.
\end{proof}

\begin{proof}[Proof of Proposition \ref{hfun_prop}]
For $\al\in\{-1/2\}\cup [1/2,\infty)$ the claim follows from \cite[Proposition~3.4]{Plewa}, hence, from now on, we consider only $\al\in(-1/2,1/2)$. Also, without any loss of generality, we assume $u\leq u'$.

Firstly, note that using the mean value theorem, Parseval's identity, and \eqref{hfun_der_estim_x>1/3} we obtain
\begin{align*}
\Vert \Rop(u,\cdot)-\Rop(u',\cdot)\Vert_{L^2(\RR_+)}&\leq |u-u'|\sup_{\xi\geq 1/2}\Vert \partial_u\Rop(\xi,\cdot)\Vert_{L^2(\RR_+)}\\
&\lesssim |u-u'|\Big( \sum_{k=0}^\infty 2^{-2k}(k+1)^{5/6}\Big)^{1/2}\\
&\lesssim |u-u'|,
\end{align*}
uniformly in $r\in(0,1/2]$ and $u,u'\geq 1/2$. On the other hand, applying \eqref{hfun_der_estim} and Lemma \ref{hfun_substr_lem}, we receive
\begin{align*}
\Big\Vert\Rop(u,\boldsymbol{\cdot})-\Rop(u',\boldsymbol{\cdot})\Big\Vert_{L^2(\RR_+)}\lesssim \sum_{k=0}^\infty 2^{-k} |\hfunk(u)-\hfunk(u')| \lesssim|u-u'|+|u-u'|^{\al+1/2},
\end{align*}
uniformly in $r\in(0,1/2]$ and $u,u'\in(0,1)$. Combining the above gives the claim for $r\in(0,1/2]$. 

Now we assume that $r\in(1/2,1)$. Invoking the formula (see \cite[p.~110]{Lebedev})
$$\frac{d}{du}I_\al(u)=\frac{\alpha}{u}I_\al(u)+I_{\al+1}(u), $$
we get
\begin{align*}
\partial_u \Rop (u,v)=\Big(\frac{2\al+1}{2u}-\frac{1+r}{1-r}u\Big)\Rop(u,v)+\frac{2rv}{1-r}R_r^{\al+1}(u,v).
\end{align*}
Using \cite[Lemma~3.2]{Plewa2} (originally from \cite[pp.~6-7]{Nasell}) we obtain
$$\Big|\frac{2rv}{1-r}R_r^{\al+1}(u,v)-\frac{1+r}{1-r}u\Rop(u,v)\Big|\lesssim \frac{1}{u}R_r^{\al+1}(u,v)+\Big(u+\frac{v-u}{1-r}\Big)\Rop(u,v), $$
uniformly in $r\in(1/2,1)$, $u,v>0$. Proceeding as in the proof of \cite[Proposition~3.4]{Plewa} one can show that
$$\Big\Vert \frac{1}{u}R_r^{\al+1}(u,\cdot)+\Big(u+\frac{\cdot-u}{1-r}\Big)\Rop(u,\cdot)\Big\Vert_{L^2(\RR_+)}\lesssim (1-r)^{-3/4}, $$
uniformly in $r\in(1/2,1)$ and $u>0$. We leave the details for the interested reader. Thus, we arrived at
\begin{align*}
\Big\Vert\Rop(u,\boldsymbol{\cdot})-\Rop(u',\boldsymbol{\cdot})\Big\Vert_{L^2(\RR_+)}&=\Big\Vert\int_{u}^{u'} \partial_s \Rop(s,\cdot)\,ds\Big\Vert_{L^2(\RR_+)}\\
&\lesssim \frac{|u-u'|}{(1-r)^{3/4}}+\Big\Vert\int_{u}^{u'}  \frac{2\al+1}{2s}\Rop(s,\cdot)\,ds\Big\Vert_{L^2(\RR_+)},
\end{align*}
uniformly in $r\in(1/2,1)$ and $u,u'>0$.

In order to complete the proof it suffices to estimate the remaining component. Using Minkowski's integral inequality and Lemma \ref{hfun_lastcomponent_lemma} we get
\begin{align*}
\Big\Vert\int_{u}^{u'}  \frac{2\al+1}{2s}\Rop(s,\cdot)\,ds\Big\Vert_{L^2(\RR_+)}&\lesssim \int_{ u}^{u'}\Big\Vert s^{-1}\Rop(s,\cdot)\Big\Vert_{L^2(\RR_+)}\,ds\\
&\lesssim |u-u'|(1-r)^{-3/4}+ (1-r)^{-(\al+1)/2}\int_{u}^{u'}s^{\al-1/2}\,ds,
\end{align*}
uniformly in $r\in(1/2,1)$ and $u,u'>0$. Finally,
\begin{equation*}
\int_{u}^{u'}s^{\al-1/2}\,ds=\int_{(u,u')\cap(0,1)}s^{\al-1/2}\,ds+\int_{(u,u')\cap(1,\infty)}s^{\al-1/2}\,ds\lesssim |u-u'|^{\al+1/2} +|u-u'|,
\end{equation*}
uniformly in $r\in(1/2,1)$ and $u,u'>0$. 

Combining the above gives the claim.
\end{proof}

{\thm\label{hfun_thm} For $\al\in[-1/2,\infty)^d$ the inequality
\begin{equation*}
\sum\limits_{n\in\NN^d}\frac{|\langle f,\hfun\rangle|}{(\ven+1)^{3d/4}}\lesssim \Vert f\Vert_{H^1(\RR_+^d)},
\end{equation*}
holds uniformly in $f\in H^1(\RR_+^d)$. The result is sharp in the sense that for any $\ve>0$ there exists $f\in H^1(\RR_+^d)$ such that
$$\sum\limits_{n\in\NN^d}\frac{|\langle f,\hfun\rangle|}{(\ven+1)^{3d/4-\ve}}=\infty. $$}
\begin{proof}
For the first part of the theorem it suffices to use Proposition \ref{hfun_prop}, \cite[Lemma~3.1]{Plewa}, and Theorem \ref{general_thm}.

In order to prove sharpness, for a given $K\in\NN$, we shall construct an appropriate $H^1(\RR^d_+)$-atom $\mathbf{a}$ such that
\begin{equation}\label{counter_atom_claim}
\sum\limits_{n\in\NN^d}^\infty\frac{|\langle \mathbf{a},\hfun\rangle|}{(\ven +1)^{3d/4-\ve}}\gtrsim K^{\ve}.
\end{equation}
We begin with the case $d=1$ and $\al>-1/2$.

Firstly, note that for $\hfunk$ we have the estimate (compare \cite[pp.~435,~453)]{Muckenhoupt})
\begin{equation*}
Ak^{\al/2}u^{\al+1/2}\leq \hfunk(u)\leq B k^{\al/2}u^{\al+1/2},\qquad 0<u\leq \frac{c}{\sqrt{k}},
\end{equation*}
where $A,B,c>0$ are constants depending only on $\al$.

Fix $\al >-1/2$, $\ve>0$, and $K\in\NN$. For $\delta\in (0,1/2)$ we define
$$a(u)=\left\{\begin{array}{ll}
\delta c^{-1}(1-\delta)^{-1} K^{1/2},& u\in(c\delta K^{-1/2},c K^{-1/2}),\\
-c^{-1}K^{1/2},& u\in(0,c\delta K^{-1/2}),
\end{array}\right. $$
It is easy to check that $a$ is an $H^1(\RR_+)$-atom. We estimate
\begin{align*}
\int_{\RR_+}a(u)\hfunk(u)\,du&\geq \frac{\delta AK^{1/2}k^{\al/2}}{c(1-\delta)}\int_{c\delta K^{-1/2}}^{cK^{-1/2}}u^{\al+1/2}\,du-\frac{K^{1/2}Bk^{\al/2}}{c} \int_{0}^{c\delta K^{-1/2}}u^{\al+1/2}\,du\\
&=\frac{2k^{\al/2}A\delta}{(2\al+3)c(1-\delta)K^{\al/2+1/4}}\left(1-\delta^{\al+1/2}(\delta+B/A)+\delta^{\al+3/2}B/A  \right)\\
&\gtrsim \frac{k^{\al/2}\delta}{K^{\al/2+1/4}(1-\delta)}\left(1- \delta^{\al+1/2}(1+B/A)\right).
\end{align*}
Choosing $\delta$ sufficiently small and independently of $K$ we obtain
$$\langle a,\hfunk\rangle\gtrsim k^{\al/2}K^{-\al/2-1/4}.$$
Thus,
\begin{align*}
\sum\limits_{k=0}^\infty\frac{|\langle a,\hfunk\rangle|}{(k+1)^{3/4-\ve}}\gtrsim K^{-\al/2-1/4}\sum\limits_{k=1}^K k^{\al/2+\ve-3/4}\gtrsim K^{\ve},
\end{align*}
which finishes the proof for $d=1$ and $\al>- 1/2$.

Note that if $\al=-1/2$, then by \eqref{derivative_formula} and \eqref{pointwise_estimates_hfun} we have
$$-\frac{d}{du}\varphi_k^{-1/2}(u) \gtrsim k^{3/4}u,\qquad 0<u\leq \frac{c}{\sqrt{k}}.$$
Hence, using the mean value theorem we obtain for $k\leq K$
\begin{align*}
\int_B a(u)\varphi_k^{-1/2}(u)\,d u&=\int_0^{c K^{-1/2}} a(u) (u-\delta cK^{-1/2}) \frac{d\varphi_k^{-1/2}}{du}(\xi_u)\,du\\
&=c^{-1}\sqrt{K}\int_0^{c K^{-1/2}}\big(\delta(1-\delta)^{-1}\ind_{B_1}(u)+\ind_{B_2}(u) \big)\big|u-\delta cK^{-1/2}\big|\\
&\qquad \times\big(-\frac{d\varphi_k^{-1/2}}{du}(\xi_u)\big)\,du\\
&\gtrsim c^{-1}\sqrt{K} k^{3/4}\delta(1-\delta)^{-1} \int_{0}^{\delta c K^{-1/2}}(\delta cK^{-1/2}-u)u\,du\\
&\simeq c^2 K^{-1}k^{3/4}\delta^3(1-\delta)^{-1}\\
&\gtrsim K^{-1}k^{3/4},
\end{align*}
where $\xi_u$ is between $u$ and $\delta c K^{-1/2}$.

In the multi-dimensional case we define
$$\mathbf{a}(x)=\prod_{i=1}^d a(x_i),\qquad x\in\RR^d_+.$$
It can be checked that $\mathbf{a}$ is an $H^1(\RR^d_+)$-atom and that \eqref{counter_atom_claim} holds.
We leave the details for the interested reader. 
\end{proof}

\section{Generalized Hermite functions}
The generalized Hermite functions of order $\lambda\geq 0$ on $\RR$ are defined by the relation
\begin{equation*}
h_{2k}^\lambda(u)=(-1)^k 2^{-1/2}\varphi_k^{\lambda-1/2}(|u|),\quad h_{2k+1}^\lambda(u)=(-1)^k2^{-1/2}\mathrm{sgn}(u)\varphi_k^{\lambda+1/2}(|u|),\quad u\in\RR,
\end{equation*}
(for $u=0$ we naturally extend the definition of $\hfunk$). In the case $d\geq 1$ we define them as tensor products of the one-dimensional $h^\lambda_k$. Note that if $\lambda=\mathbf{0}$, then the functions $\{h^{\mathbf{0}}_n\}_{n\in\NN^d}$ are the classical Hermite functions.

In the following theorem we use two inner products: in $L^2(\RR^d)$ and in $L^2(\RR^d_+)$  denoted by $\langle \cdot,\cdot\rangle$ and $\langle \cdot,\cdot\rangle_+$, respectively.

{\thm\label{Hfun_H1_thm} Let $\lambda\in [0,\infty)^d$. The following inequality holds
$$\sum_{n\in\NN^d}\frac{|\langle f,h_n^{\lambda}\rangle|}{(\ven+1)^{3d/4}}\lesssim \Vert f\Vert_{H^1(\RR^d)}, $$
uniformly in $f\in H^1(\RR^d)$. The exponent is sharp, in the sense that for every $\ve>0$ there exists $f\in H^1(\RR^d)$ such that
$$\sum_{n\in\NN^d}\frac{|\langle f,h_n^{\lambda}\rangle|}{(\ven+1)^{3d/4-\ve}}=\infty. $$}
\begin{proof}
We shall justify that the claims follow from Theorem \ref{hfun_thm}.

We introduce a function $\mathfrak{m}\colon \NN^d\rightarrow \{0,1\}^d$, defined by 
$$\mathfrak{m}(n)_i=n_i\quad(\mathrm{mod}\,2),\qquad i=1,\ldots,d.$$
Fix $\lambda\in[0,\infty)^d$. For $\eta\in\{0,1\}^d$ we shall denote
$$\lambda_\eta=\Big(\lambda_1-\frac{(-1)^{\eta_1}}{2},\ldots,\lambda_d-\frac{(-1)^{\eta_d}}{2}\Big). $$
Note that $h^\lambda_n$ is $\mf(n)$-symmetric. Hence, 
$$\langle f,h_n^\lambda\rangle\simeq \big\langle f_{\mf(n)}^+,\varphi_{\lfloor n/2\rfloor}^{\lambda_{\mf(n)}}\big\rangle_+,\qquad f\in H^1(\RR^d),\ n\in\NN^d. $$
Thus, we estimate using \eqref{H1_norm_estimate}
\begin{align*}
\sum_{n\in\NN^d}\frac{|\langle f,h_n^{\lambda}\rangle|}{(\ven+1)^{3d/4}}\simeq\sum_{\eta\in\{0,1\}^d} \sum_{\mf(n)=\eta}\frac{|\langle f_\eta^+,\varphi_{\lfloor n/2\rfloor}^{\lambda_\eta}\rangle_+|}{(\ven+1)^{3d/4}}&\leq \sum_{\eta\in\{0,1\}^d} \sum_{n\in\NN^d}\frac{|\langle f_\eta^+,\varphi_n^{\lambda_\eta}\rangle_+|}{(\ven+1)^{3d/4}}\\
&\lesssim \sum_{\eta\in\{0,1\}^d} \Vert f_\eta^+\Vert_{H^1(\RR^d_+)}\\
&\lesssim \Vert f\Vert_{H^1(\RR^d)}.
\end{align*}
This finishes the verification of the first claim.

In order to prove the second claim, we fix $\ve>0$. Let $\al=\lambda-\mathbf{1/2}$. Theorem \ref{hfun_thm} yields that there exists  $g\in H^1(\RR^d_+)$ such that
\begin{equation*}
\sum\limits_{n\in\NN^d}\frac{|\langle g,\hfun\rangle_+|}{(\ven+1)^{3d/4-\ve}}=\infty.
\end{equation*}
We extend $g$ to an $\mathbf{0}$-symmetric function $f$. We emphasise that $f\in H^1(\RR^d)$. Hence,
\begin{align*}
\sum_{n\in\NN^d}\frac{|\langle f,h_n^{\lambda}\rangle|}{(\ven+1)^{3d/4-\ve}}\geq \sum_{\mathfrak{m}(n)=\mathbf{0}}\frac{|\langle f,h_n^{\lambda}\rangle|}{(\ven+1)^{3d/4-\ve}}
\simeq \sum_{n\in\NN^d}\frac{|\langle g,\varphi_{n}^{\lambda-\mathbf{1/2}}\rangle_+|}{(2^d\ven+1)^{3d/4-\ve}}=\infty.
\end{align*}
This finishes the proof of the theorem.
\end{proof} 

Theorem \ref{Hfun_H1_thm} holds for the classical Hermite functions (that is for $\lambda=\mathbf{0}$), and hence the admissible exponent obtained in \cite{LiYuShi, RadhaThangavelu} is sharp. 

In the previous articles (see \cite{Plewa,Plewa2}) we proved the $L^1$-analogues of Hardy's type inequalities. Therefore we present a corresponding result for the generalized Hermite functions below. It can be proved basing on \cite[Theorem~5.1]{Plewa} and using similar arguments as in the proof of Theorem \ref{Hfun_H1_thm}.

{\thm\label{Hfun_L1_thm} Let $\ve>0$ and $\lambda\in [0,\infty)^d$. Then 
\begin{equation*}
\sum_{n\in\NN^d}\frac{\vert\langle f,h_n^\lambda\rangle\vert}{(\ven+1)^{3d/4+\ve}}\lesssim \Vert f\Vert_{L^1(\RR^d)},
\end{equation*}
uniformly in $f\in L^1(\RR^d)$.
The result is sharp in the sense that there is $f\in L^1(\RR^d)$ such that
\begin{equation*}
\sum_{n\in\NN^d}\frac{\vert\langle f,h_n^\lambda\rangle\vert}{(\vert n\vert+1)^{3d/4}}=\infty.
\end{equation*}
}

\end{document}